\documentclass[12pt,a4paper]{article}
\usepackage{mathrsfs}
\usepackage{amsmath,amssymb}
\usepackage{latexsym}
\usepackage{a4wide}
\usepackage{graphicx}
\usepackage{subfigure}
\usepackage{caption}

\begin{document}
\graphicspath{{Latex_Graphics/}}

\title{{\Large A New Boundary Scheme for BGK Model}\thanks{This
work is partially supported by the Natural Science Foundation of
China(11071175).}}

\author{ {\normalsize Peng Wang\footnote{\quad{\small {email: 7790326@qq.com. }}}, Shiqing Zhang\footnote{\quad{\small
email: Zhangshiqing@msn.com. }}}\\
{\small Yangtze Center of Mathematics and College of Mathematics,}\\
{\small Sichuan University, Chengdu 610064, People's Republic of
China} }
\date{}
\maketitle

\begin{abstract}
In this paper, we proposed a new boundary condition scheme for the
BGK model which has second order accuracy and can be developed to
higher accuracy. Numerical tests show that the numerical solutions
of the BGK model applied to the new boundary scheme have great
agreement with analytical solutions. It is also found that the
numerical accuracy of the present schemes is much better than
that of the original extrapolation schemes proposed by Guo et al at
small grid.

\end{abstract}
%\tableofcontents
%\parskip 1pt
\section{Introduction}

\hspace*{\parindent}In the past 20 years, the Boltzmann-BGK method
in the simulation of fluid showed great success.  Compared to
traditional numerical methods, the Boltzmann-BGK method that based
on the mesoscopic view got macroscopic variables throughout the
moment integration of distribution function $f$.  For the
mesoscopic property, the Boltzmann-BGK can be easily applied to a very
large scope, such as multiphase, multicomponent high-speed
compressible flows and so on\cite{BOOKHe1,BOOKGuo1,Qi2001107}. In simulation, the boundary condition in practice is usually given
by macroscopic variables, in BGK method we actually use the
distribution function. Up to now there are no universal method to
get the distribution function throughout the macroscopic variables,
so it's necessary to discuss the distribution function of
boundary conditions.

Virtual equilibrium method, interpolation method  and
non-equilibrium extrapolation method are the main boundary schemes
used in BGK model nowadays. Most of the schemes based on some sort
of extrapolation methods to get the distribution functions on boundary
\cite{Filippova1998219,1009-1963-11-4-310,Mei1998426,Mei1999307,bouzidi:3452}. In this paper, we got
a new boundary scheme based on mathematical analysis. This new
boundary scheme has second order accuracy, can be easily to develop
to higher order accuracy and can be implemented by different viscosity
between flow node and the boundary node. In the following sections, we
first give implicit-explicit scheme of BGK model, and then discuss
boundary scheme, in the last, we applied this new boundary condition
scheme together with the implicit-explicit scheme of BGK model to
Couette flow and lid driven cavity flow to test accuracy and
stability of this boundary scheme.

\section{BGK method}
\subsection{Basic BGK model}

\hspace*{\parindent}The standard dynamics theory of mesoscopic model
is described by the Boltzmann equation as follow
\begin{equation}\label{Boltzmann_equation}
\frac{\partial f}{\partial t} + \vec{e} \cdot \nabla f = J,
\end{equation}
where $f$ denotes the distribution function, $\vec{e}$
denotes the molecule velocity vector, and $J$ represents the
collision term\cite{cercignani1988boltzmann}. By
the model of BGK, we can easily simplify the collision term $J$ to
the following formulation
\begin{equation}\label{continuous_bgk_equation}
\frac{\partial f}{\partial t} + \vec{e} \cdot \nabla f =
\frac{1}{\tau}(f^{(eq)} - f  ),
\end{equation}
where $\tau$ denotes the relaxation time and $f^{(eq)}$ denotes the
local equilibrium distribution function\cite{PhysRev.94.511}.

By the Hermite expansion or Taylor expansion, we discrete the
infinite velocity space to finite space $\{e_0,e_2,\cdots,e_N\}$ to
make the equation (\ref{continuous_bgk_equation}) be a set of
discrete velocity equations.
\begin{equation}\label{BGK_dve}
\frac{\partial f_i}{\partial t} + \vec{e_i} \cdot \nabla f_i =
\frac{1}{\tau}(f_i^{(eq)} - f_i),(i=0,1,\cdots,N).
\end{equation}
The macroscopic density $\rho$ and velocity $\vec{u}$ can be
obtained by the moments of the distribution function
\begin{equation}\label{macroscopic_variavles}
\rho = \sum f_i ,\,\,\,\, \rho \vec{u} = \sum \vec{e}_i f_i.
\end{equation}

 In 1992, Qian proposed a D2Q9 model\cite{0295-5075-17-6-001}. The $f_i^{(eq)}$ of D2Q9
 is
\begin{equation}\label{feq_form}
f_i^{(eq)} = \rho \omega_i [1 + \frac{\vec{e}_i\cdot \vec{u}}{c_s^2}
+ \frac{(\vec{e}_i\cdot \vec{u})^2}{2c_s^4} -
\frac{|\vec{u}|^2}{2c_s^2}].%\,\,\,\,(i = 0,1,...,8)
\end{equation}
where $\omega_i$ is a weight coefficient, $c_s = 1/\sqrt{3}$ is the
local sound speed. For D2Q9 model, the discrete speed is
\begin{equation}
\vec{e}_i = c\begin{pmatrix}
            0 & 1 & 0 & -1 & 0 & \sqrt{2} & -\sqrt{2} & -\sqrt{2} & \sqrt{2} \\
            0 & 0 & 1 & 0 & -1 & \sqrt{2} & \sqrt{2} & -\sqrt{2} & -\sqrt{2} \\
          \end{pmatrix}\,\,\,\,(i = 0,1,...,8),
\end{equation}
where the $c=1$. The sound speed and weight coefficient is
\begin{equation}
c_s = \frac{1}{\sqrt{3}}, \omega_i = \left\{ \begin{aligned}
                                            4/9,\,\,\,\,\,\,\,\, &|\vec{e}_i|^2= 0, \\
                                            1/9,\,\,\,\,\,\,\,\, &|\vec{e}_i|^2=c^2, \\
                                            1/36,\,\,\,\,\,\,\,\,&|\vec{e}_i|^2=2c^2.\\
                                            \end{aligned} \right.
\end{equation}

Apply Enskog-Chapman expansion to recover the macroscopic
incompressible Navier-Stokes equations without outer force
\begin{equation}
\begin{aligned}
& \nabla \cdot \vec{u} = 0,\\
&\frac{\partial\vec{u}}{\partial t} +\vec{u}\cdot\nabla \vec{u} =
-\frac{1}{\rho}\nabla p + \mu \Delta \vec{u},
\end{aligned}
\end{equation}
where the $\mu$ represents the kinematic viscosity. We need the
discrete velocity and equilibrium distribution function satisfying the
following equation:
\begin{equation}
\rho = \sum_{i}{f^{(eq)}_i},\rho \vec{u} = \sum_{i}{\vec{e}_i
f^{(eq)}_i},\rho u_i u_j + p\delta_{ij} =
\sum_{i}{e_{i\alpha}e_{i\beta} f^{(eq)}_i},
\end{equation}
and
\begin{equation}
\mu = \tau c_s^2.
\end{equation}

In this paper, we applied an explicit-implicit scheme to test the
validness of the new scheme of boundary condition.

\subsection{Time and spatial discrete} \hspace*{\parindent}For the
calculation of the equations (\ref{BGK_dve}) of the basic BGK model,
we integrate the equation (\ref{BGK_dve}) in both side to get
\begin{equation}\label{t1}
\left\{
\begin{aligned}
    &\int_{t}^{t + \Delta t}{\frac{\partial f_i^{}}{\partial t}}+ \int_{t}^{t + \Delta
    t}{\vec{e}_i\cdot \nabla f_i} = \int_{t}^{t + \Delta t}{\frac{1}{\tau}(f_i^{(eq)} -
    f_i)},\\
    &f(0,\vec{x},\vec{e}) = f_0(\vec{x},\vec{e}).
\end{aligned} \right.
\end{equation}
where the $f_0(\vec{x},\vec{e})$ represents the initial value. By
mean value theorem of integrals, the equality (\ref{t1}) can be
written as
\begin{equation}\label{mean_theroy}
f_i^{n+1} - f_i^{n} + \vec{e}_i\cdot \nabla f_i(t_\xi,\vec{x})
\Delta t = \frac{\Delta t}{\tau}(f_i^{(eq)}(t_\xi,\vec{x}) -
    f_i(t_\xi,\vec{x})),
\end{equation}
where $t<t_\xi<t+\Delta t$, the superscript $n$ represents the step
number, i.e. $f_i^n = f(t,\vec{x}),f_i^{n+1} = f(t+\Delta
t,\vec{x})$.

For time discrete, we calculates advection term $\vec{e}_i\cdot
\nabla f_i(t_\xi,\vec{x}) \Delta t$ as an explicit finite-difference
form and the collision term $ \frac{\Delta
t}{\tau}(f_i^{(eq)}(t_\xi,\vec{x}) -
    f_i(t_\xi,\vec{x}))$ as a
implicit finite-difference form, and we got
\begin{equation}\label{iterative_scheme}
f_i^{n+1} - f_i^n + \Delta t \vec{e}_i \cdot \nabla f_i^n = \Delta t
[\theta J_i^{n+1} + (1-\theta)J_i^n],
\end{equation}
where the $\theta$ represents the degree of implicity, in this paper
we make it to be $0.5$, and $J_i^n = \frac{\Delta
t}{\tau}(f_i^{(eq)}(t,\vec{x}) -
    f_i(t,\vec{x})),\, J_i^{n+1} = \frac{\Delta
t}{\tau}(f_i^{(eq)}(t+\Delta t,\vec{x}) -
    f_i(t+\Delta t,\vec{x})).$

Z.L Guo and T.S Zhao introduced a new distribution function to remove
the implicity of the equation (\ref{iterative_scheme}), and the new distribution function is
\begin{equation}\label{new_df}
g_i = f_i +   \pi\theta(f_i - f_i^{(eq)})
\end{equation}
where $\pi = \Delta t/\tau$\cite{PhysRevE.67.066709}. Applying this new distribution function
to equation (\ref{iterative_scheme}), we can delete the implicity
of the collision term and got
\begin{equation}
\begin{aligned}
&f_i^n = \frac{1}{1+ \pi \theta} (g_i^n + \pi \theta f^{(eq),n}),\\
&g_i^{n+1} = -\Delta t \vec{e}_i\cdot \nabla f_i^n + (1-\pi + \pi
\theta)f_i^{n} + \pi(1-\theta)f_i^{(eq),n}.
\end{aligned}
\end{equation}

In the following part, we discuss about the spatial discrete scheme.
We apply a mixed-difference scheme which combined upwind scheme
with central scheme. The central difference of $\frac{\partial f_i}{\partial
x_\alpha}$ is
\begin{equation}
(\frac{\partial f_i}{\partial x_\alpha} )_c= \frac{1}{2\Delta
x_\alpha} [f_i(x_\alpha + \Delta x_\alpha, \cdot) - f_i(x_\alpha -
\Delta x_\alpha, \cdot)].
\end{equation}
The second-order upwind-difference scheme is
\begin{equation}
(\frac{\partial f_i}{\partial x_\alpha})_u = \left \{
\begin{aligned}
&\frac{1}{2\Delta x_\alpha} [3f_i(x_\alpha,\cdot) - 4f_i(x_\alpha -
\Delta
x_\alpha, \cdot) + f_i(x_\alpha - 2\Delta x_\alpha, \cdot)] \,\,\, if\, e_{i\alpha}\geqslant 0,\\
&-\frac{1}{2\Delta x_\alpha} [3f_i(x_\alpha,\cdot) - 4f_i(x_\alpha +
\Delta x_\alpha, \cdot) + f_i(x_\alpha + 2 \Delta x_\alpha,
\cdot)]\,\,\, if\,e_{i\alpha} < 0 .
\end{aligned}
\right.
\end{equation}
The mixture form is
\begin{equation}\label{space_mitxture_form}
\frac{\partial f_i}{\partial x_\alpha} = (\zeta\frac{\partial
f_i}{\partial x_\alpha})_c + (1-\zeta)(\frac{\partial f_i}{\partial
x_\alpha})_u,
\end{equation}
where $\zeta \in [0,1]$.

By the time discrete and the space discrete, we got the entire
iterative form:
\begin{equation}\label{entire_iterative_form}
\begin{aligned}
&g_i^{n+1} = -\Delta t \vec{e}_i\cdot \nabla_h f_i^n +
(1-\pi + \pi \theta)f_i^{n} + \pi(1-\theta)f_i^{(eq),n},\\
&f_i^n = \frac{1}{1+ \pi \theta} (g_i^n + \pi \theta f^{(eq),n}),
\end{aligned}
\end{equation}
where the $\nabla_h$ represents the euqtion
(\ref{space_mitxture_form}). By equation (\ref{new_df}) and
(\ref{macroscopic_variavles}), we got
\begin{equation}
\rho^{n+1} = \sum g_i^{n+1},\,\,\,\, \vec{u}^{n+1} =
\frac{1}{\rho^{n+1}}\sum \vec{e}_i g_i^{n+1}.
\end{equation}

\subsection{Scheme of boundary condition}
\hspace*{\parindent}It's well-known that the scheme of the boundary
condition plays a very important role in the flow computation. It not
only relates to the stability of the computation, but also relates
to the the accuracy of the computation. In simulation, the boundary
condition is usually given by macroscopic variables, but in BGK
method, we actually use the distribution function. For those reasons,
we discuss the following boundary conditional scheme for the BGK model. Be
Inspired by non-equilibrium extrapolation method which was proposed
by Guo, Zheng and Shi\cite{1009-1963-11-4-310}, we designed the new scheme of boundary condition. The basic
idea of Guo et al is decomposed the distribution function on the
boundary into equilibrium part and the non-equilibrium part
\begin{equation}
f_i(t,\vec{x}_p) = f_i^{(eq)}(t,\vec{x}_p) +
f_i^{(neq)}(t,\vec{x}_p).
\end{equation}
where the $\vec{x}_p$ denotes the position of the physic boundary. At present work, we also decompose boundary distribution functions
into equilibrium part and the non-equilibrium part. In BGK model, the
equilibrium $f_i^{(eq)}(t,\vec{x}_p)$ functions can be assumed as the
functional of $\vec{u},\rho,p,T$ which are related with the time and space\cite{cercignani1988boltzmann,BOOKLi2}. For this reason
$f_i^{(eq)}(t,\vec{x}_p)$ was approximated by
\begin{equation}\label{feq_approximationg}
f_i^{(eq)}(t,\vec{x}_p)=
f_i^{(eq)}(t,\rho(t,\vec{x}_f),\vec{u}(t,\vec{x}_p)).
\end{equation}
where the $\vec{x}_f$ denotes the flow node that next to the
$\vec{x}_p$ i.e $\vec{x}_f = \vec{x}_p + \vec{e}_i\Delta x$, and
this approximation is at least third order
accuracy\cite{1009-1963-11-4-310}. For
the $f_i^{(neq)}(t,\vec{x}_p)$ part, consider the boundary
distribution function which satisfies the equation
\begin{equation}\label{BGK_dve_boundary}
\frac{\partial f_i(t,\vec{x}_p)}{\partial t} + \vec{e_i} \cdot
\nabla f_i(t,\vec{x}_p) = \frac{1}{\tau_*}(f_i^{(eq)}(t,\vec{x}_p) -
f_i(t,\vec{x}_p)),(i=0,1,\cdots,N).
\end{equation}
where the $\tau_*$ is determined by the viscosity between physical
boundary and the fluid. If we think that the viscosity between
physical boundary and the fluid is the same as inner flow, we have
the $\tau_*=\tau$. In this paper, we set $\tau_* = \tau$.
Assume the solution of equation (\ref{BGK_dve_boundary}) can be
expanded as
\begin{equation}\label{boundary_solution_expand}
f_i(t,\vec{x}_p) = f_i^{(eq)}(t,\vec{x}_p) + \tau_*
f_i^{(1)}(t,\vec{x}_p) + \tau^2_* f_i^{(2)}(t,\vec{x}_p) + \cdots.
\end{equation}
Substitute (\ref{boundary_solution_expand}) into equation
(\ref{BGK_dve_boundary}) and compare the for coefficient
$\tau_*$, we have
\begin{equation}\label{f_1_form}
f_i^{(1)}(t,\vec{x}_p) = -(\frac{\partial }{\partial t} +
\vec{e}_i\cdot \nabla )f_i^{(eq)}(t,\vec{x}_p).
\end{equation}

Based on equation (\ref{feq_form}) and (\ref{feq_approximationg}), we
can easily got the $f_i^{(eq)}(t,\vec{x}_p)$ throughout the
macroscopic variables. By $f_i^{(eq)}(t,\vec{x}_p)$, we
can got the $f^{(1)}(t,\vec{x}_p)$ based on equation
(\ref{f_1_form}). We use $\tau_* f_i^{(1)}(t,\vec{x}_p)$ to
approximate $f^{(neq)}(t,\vec{x}_p)$ part. By equation
(\ref{boundary_solution_expand}) , the accuracy of
$f_i^{(eq)}(t,\vec{x}_p) + \tau_* f_i^{(1)}(t,\vec{x}_p)$ to
approximate $f_i(t,\vec{x}_p)$ is second order. On the boundary, we have distribution
\begin{equation}\label{Boundary_scheme}
f_i(t,\vec{x}_p) = f_i^{(eq)}(t,\vec{x}_p) - \tau_* (\frac{\partial
}{\partial t} + \vec{e}_i\cdot \nabla )f_i^{(eq)}(t,\vec{x}_p).
\end{equation}
%where the $f_i^{(eq)}(t,\vec{x}_p)=
%f_i^{(eq)}(t,\rho(t,\vec{x}_f),\vec{u}(t,\vec{x}_p))$.

If we want to get the higher accuracy scheme, we can mix Guo's scheme with
proposed scheme, and we have
\begin{equation}\label{higher_Boundary_scheme}
\begin{aligned}
f_i(t,\vec{x}_p) =& f_i^{(eq)}(t,\vec{x}_p) - \tau_* (\frac{\partial
}{\partial t} + \vec{e}_i\cdot \nabla )f_i^{(eq)}(t,\vec{x}_p)\\
 &+
f_i(t,\vec{x}_f) + \tau_* (\frac{\partial }{\partial t} +
\vec{e}_i\cdot \nabla )f_i^{(eq)}(t,\vec{x}_f) -
f_i^{(eq)}(t,\vec{x}_f),
\end{aligned}
\end{equation}
%or we have
%\begin{equation}\label{higher_Boundary_scheme_w}
%\begin{aligned}
%f_i(t,\vec{x}_p) =& f_i^{(eq)}(t,\vec{x}_p) - \tau_* (\frac{\partial
%}{\partial t} + \vec{e}_i\cdot \nabla )f_i^{(eq)}(t,\vec{x}_p)\\
% &+
%f_i(t,\vec{x}_f) + \tau_* (\frac{\partial }{\partial t} +
%\vec{e}_i\cdot \nabla )f_i^{(eq)}(t,\vec{x}_f) -
%f_i^{(eq)}(t,\vec{x}_f).
%\end{aligned}
%\end{equation}
To get the value of $f^{(1)}(t,\vec{x}_p)$, we discrete the time
differential operator $\frac{\partial }{\partial t}$ as
\begin{equation}\label{Boundary_scheme_time}
(\frac{\partial f_i^{(eq)}(t,\vec{x}_p)}{\partial t} )_{BT} =
\frac{1}{\Delta t} [ f_i^{(eq)}(t+ \Delta t,\cdot) - f_i^{(eq)}(t,
\cdot)],
\end{equation}
and discrete differential operator $\nabla$ at boundary as
%\begin{equation}\label{Boundary_scheme}
%(\frac{\partial f_i^{(eq)}}{\partial x_\alpha} )_{Boundary} =
%\frac{1}{\Delta x_\alpha} [ f_i^{(eq)}(t,x_\alpha, \cdot) -
% f_i^{(eq)}(t,x_\alpha - \Delta x_\alpha, \cdot)],
%\end{equation}
\begin{equation}\label{Boundary_scheme_BS}
(\frac{\partial f_i^{(eq)}(t,\vec{x}_p)}{\partial x_\alpha} )_{BS} =
\frac{1}{2\Delta x_\alpha} [ 3f_i^{(eq)}(x_\alpha, \cdot) -
 4 f_i^{(eq)}(x_\alpha - \Delta x_\alpha, \cdot) + f_i^{(eq)}(x_\alpha - 2\Delta
x_\alpha, \cdot)].
\end{equation}

Summing up the above discussion, we have entire finite difference
method for the Blotzmann-BGK model. The following part, we will apply above
boundary scheme together with the implicit-explicit
difference scheme to test the validness of the proposed boundary
scheme.

\section{Numerical Simulations}
\hspace*{\parindent}In this section, we will apply the new
 boundary scheme together with the implicit-explicit scheme of BGK model to Couette flow and lid driven
cavity flow to test accuracy and stability of this boundary scheme.

\subsection{Couette flow}
\hspace*{\parindent}  The Couette plate flow is defined in the
region $0\leqslant x \leqslant 1,0\leqslant y \leqslant 1$ under a
periodic boundary condition at the entrance and exit. The bottom
plate is kept stationary, and the top plate moves horizontally with
a constant velocity $u_0$. This Couette flow has the following
analytical solution
\begin{equation}
\vec{u}^*(t,x,y) = (\frac{y}{H} + 2\sum_{k = 1}^{\infty}
\frac{(-1)^k}{\lambda_k H} \exp{(-\mu \lambda_k^{2}
t)}\sin{(\lambda_k y),0)}
\end{equation}
where the $\lambda_k = k \pi/H$, $k = 1,2...$ and $H = L = 1.0$ . We
defined the average error as
\begin{equation}\label{Average_error}
AverageError = \frac{1}{n}\sum_{n}{\frac{\sqrt{(u_{x,n} -
u^*_{x,n})^2 + (u_{y,n} - u^*_{y,n})^2}}{\sqrt{( u^*_{x,n})^2 +
(u^*_{y,n})^2}}}
\end{equation}
where the $n$ denotes the number of the grid.

\begin{figure}[!h]
\begin{minipage}[t]{0.5\linewidth}
\centering
\includegraphics[scale=0.35,trim = 0 0 0 85mm]{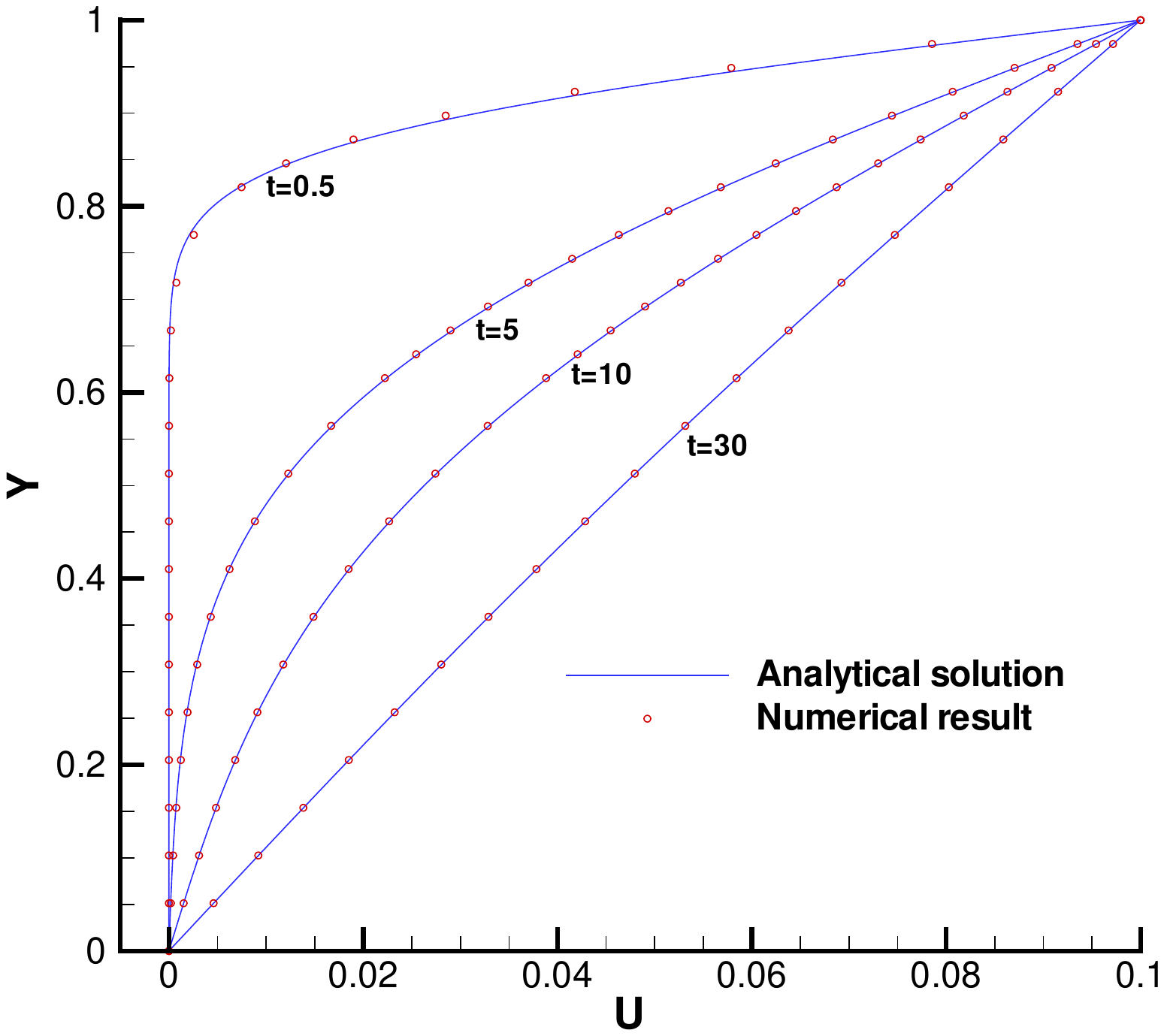}
\caption{Numerical result which applied the proposed boundary scheme
for the couette flow} \label{figure_profile_couette}
\end{minipage}
\begin{minipage}[t]{0.5\linewidth}
\centering
\includegraphics[scale=0.35,trim = 0 0 0 85mm]{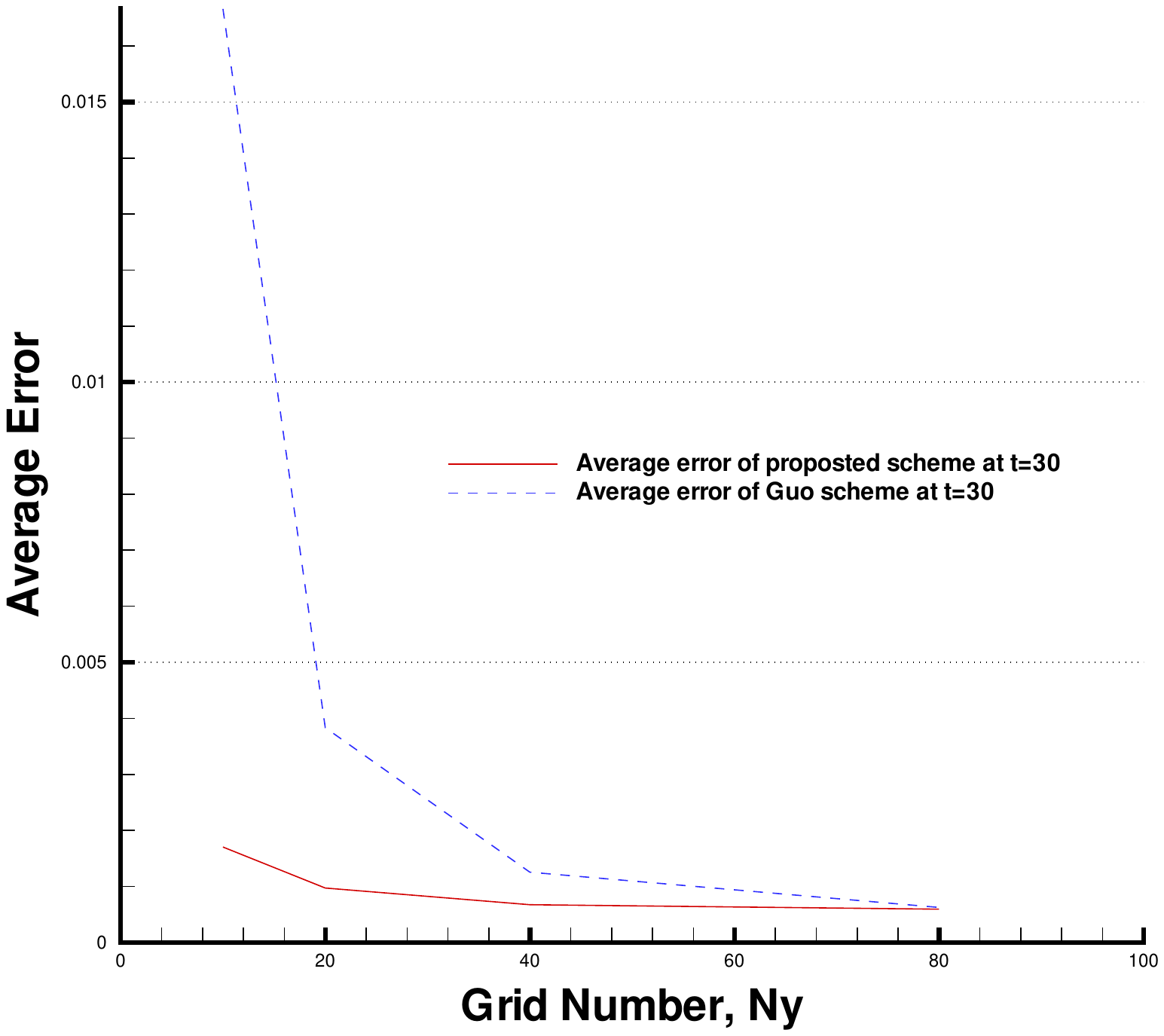}
\caption{Average error of Guo scheme and proposed scheme for couette
flow } \label{Couette_average_error}
\end{minipage}
\end{figure}

 We set $Re= (Lu_0)/\mu = 10$, $\zeta = 0.9$ and use the grid $N_x \times N_y =5 \times 10,N_x \times N_y =10 \times 20,N_x \times N_y =20 \times 40,N_x \times N_y =40 \times
 80$ for simulation. The proposed scheme (\ref{Boundary_scheme}) is
 implemented to top and bottom plates, and the periodic boundary
 condition is implemented to the entrance and
 exit\cite{BOOKHe1}. The numerical velocity profiles at time $t= 0.5,t=5,t=10,t=30$ together with analytical solutions are plotted in
 figure \ref{figure_profile_couette}. From the figure, we can see the
 numerical solution greatly agreed with the analytical solution. The
 average error of Guo scheme and proposed scheme are plotted in
 figure \ref{Couette_average_error}. Comparing Guo's scheme, the
 proposed scheme at small grid have better accuracy.

\subsection{Two dimensional lid driven square
cavity flow}

\hspace*{\parindent}The lid driven square cavity flow is classical
benchmark problem in numerical simulation of the fluid. In lid
driven square cavity flow, the top lid of the cavity has a constant
velocity toward right, and the other three boundary hold still. The
geometry of the lid driven flow is very simply and the phenomena is
very complicated.

\begin{figure}[!h]
\begin{center}
\includegraphics[scale=0.35,trim = 0 0 0 140mm]{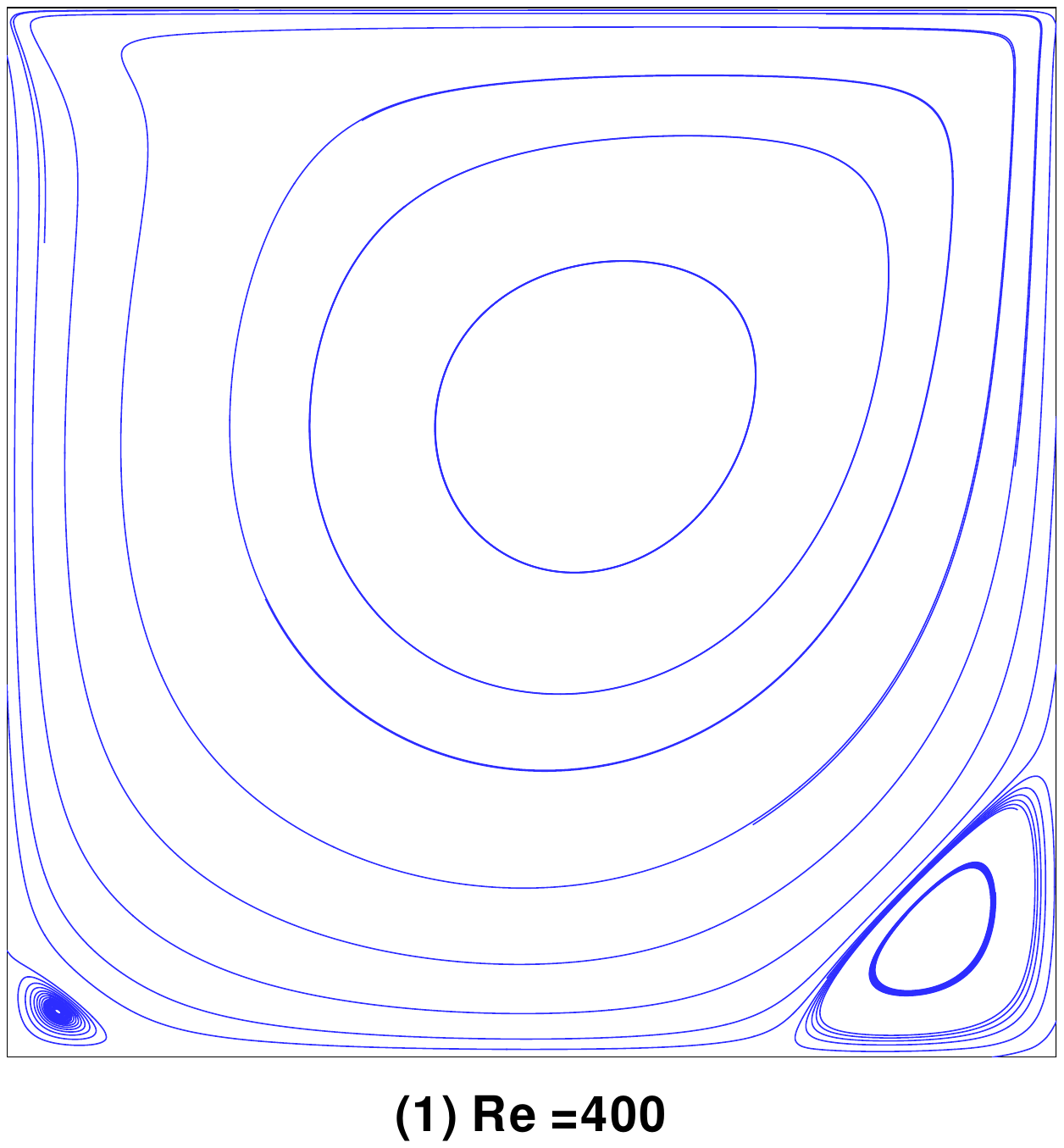}\includegraphics[scale=0.35,trim = 0 0 0 140mm]{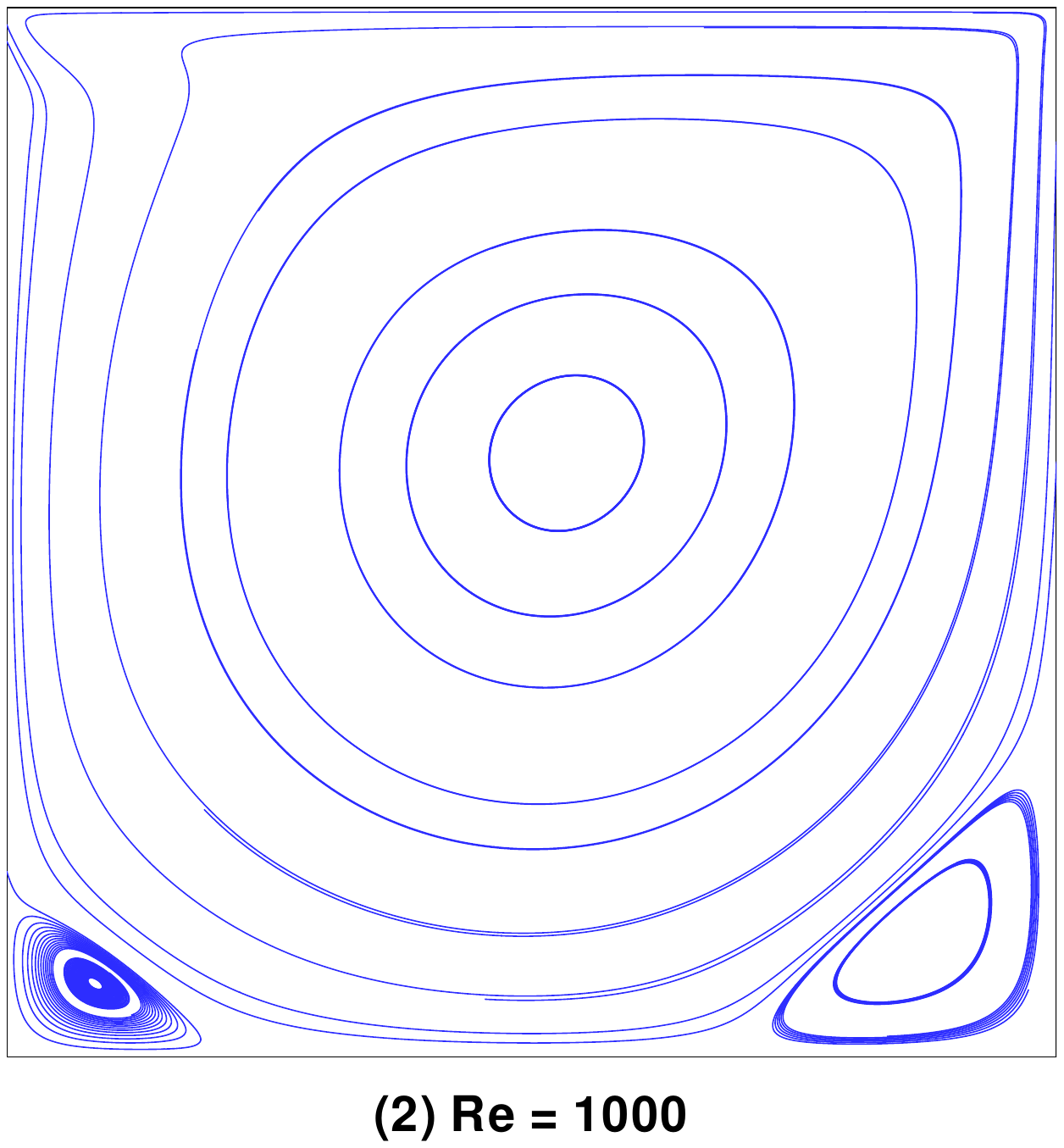}
\includegraphics[scale=0.35,trim = 0 0 0 140mm]{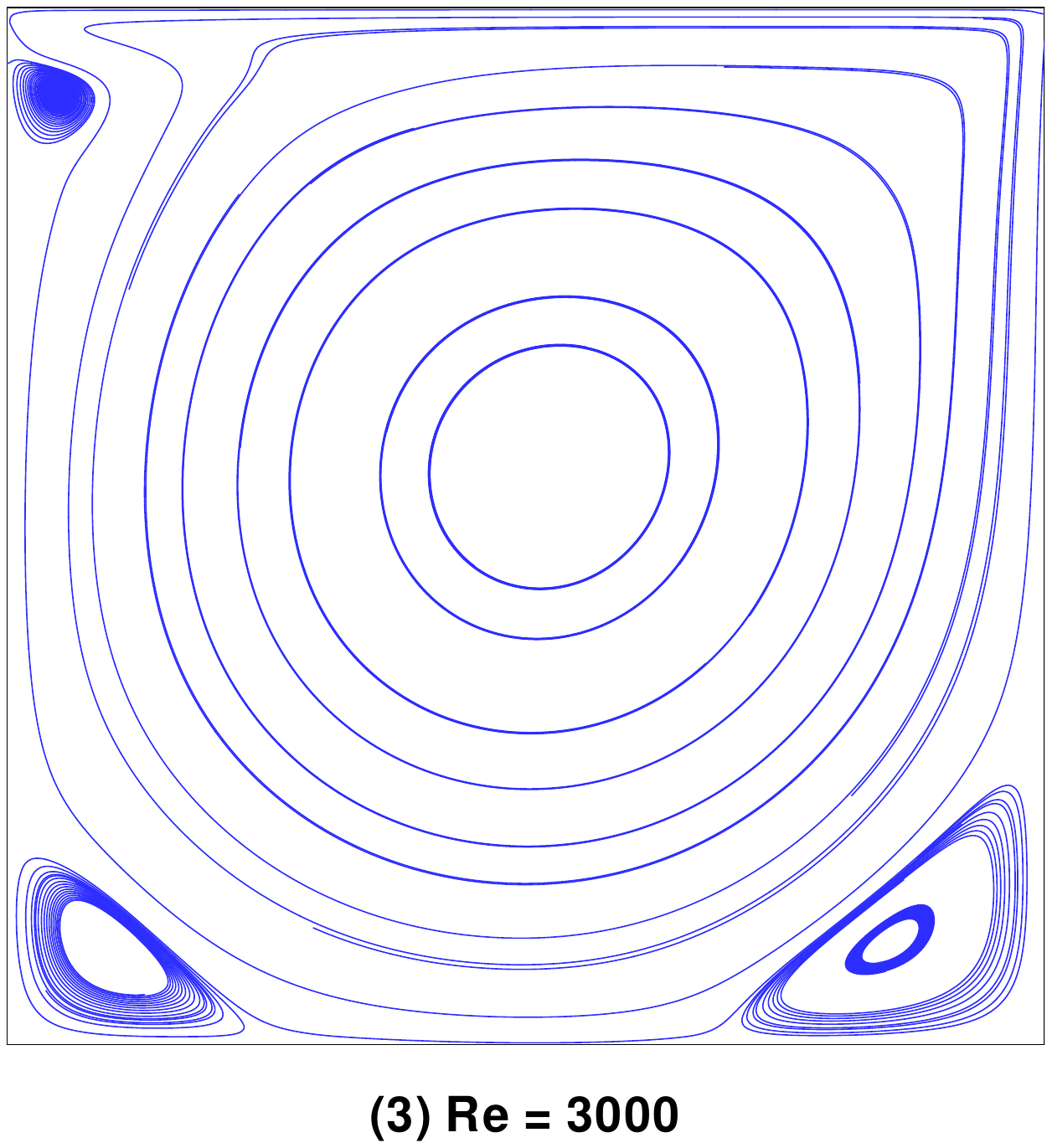}\includegraphics[scale=0.35,trim = 0 0 0 140mm]{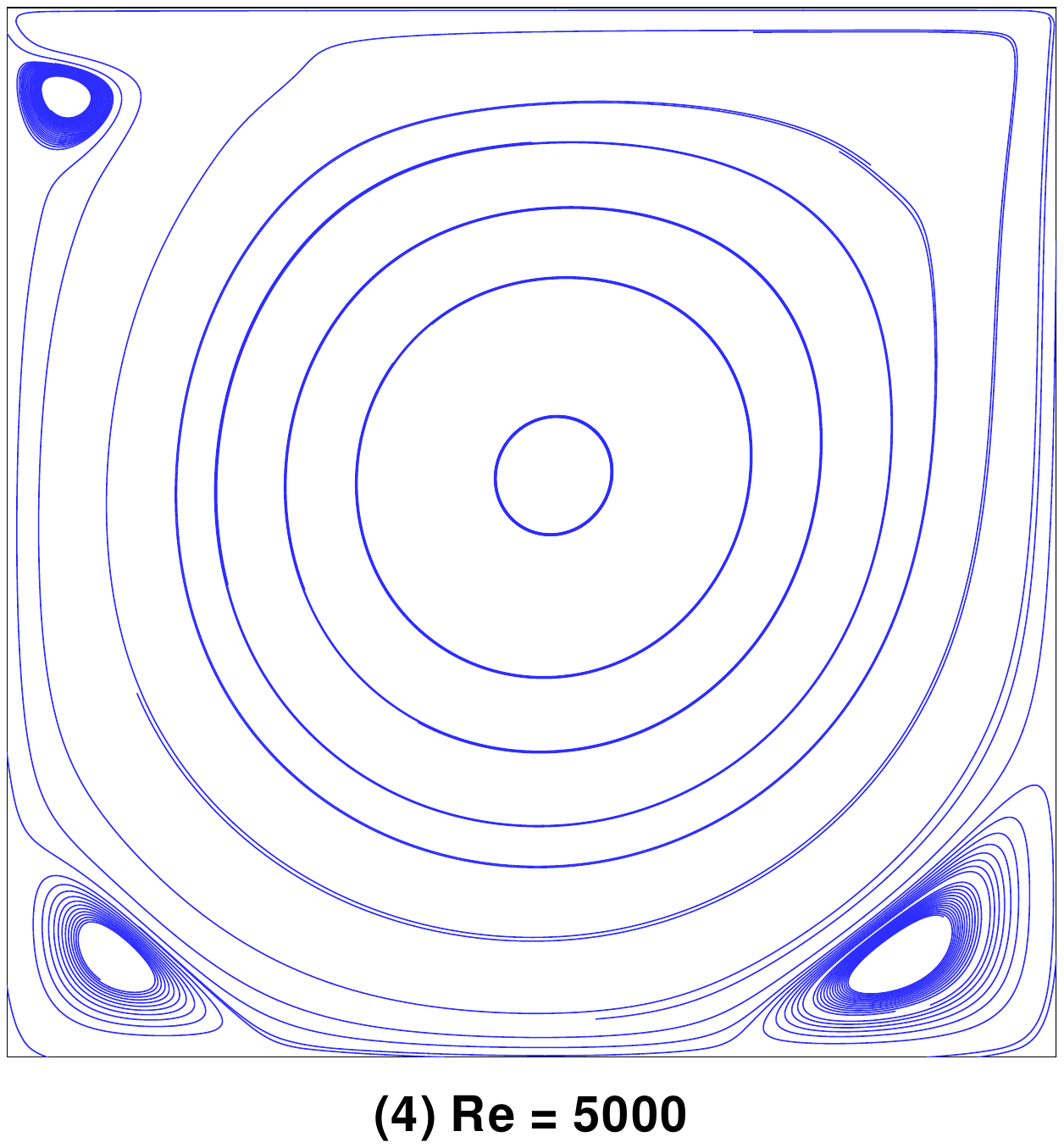}
\caption{Streamline of the lid driven flow at different Reynolds
numbers.
(1)Re=400;(2)Re=1000;(3)Re=3000;(4)Re=5000.}\label{figure_lid}
\end{center}
\end{figure}

In this section we apply proposed boundary scheme for the lid
driven cavity flow. The width and height is chosen to be $L=1.0$.
The initial density $\rho = 1.0$, top velocity of lid is chosen to
be $\vec{u}=(u_x,u_y)=(0.1,0)$. The size of the mesh is $N_x \times
N_y = 128 \times 128 $, set $\Delta t = 0.1 \times \Delta y$.
The streamline for $Re = 400, 1000, 3000, 5000$ were plotted in
figure \ref{figure_lid}.
For $Re = 400$, we can see three vortices in the figure and the center
of the vortex is at $x= 0.5563, y =0.6103$. For $Re = 1000$ it also
has three vortices and the coordinates of the vortex center is $x=0.5334,
y =0.5759$. For $Re = 3000$, we see four vortices in the picture and the
center of the vortex is $x=0.5224, y =0.5581$. When  $Re=5000$, we can see there are four vortices in the figure and the center of
the vortex is $x=0.5210,y=0.5543$. All four position of the centers of
the vortex well matched the paper\cite{BOOKHe1,PhysRevE.67.066709}.

\section{Conclusion}
\hspace*{\parindent}In this paper, we discussed a new boundary
scheme together with the implicit-explicit scheme for BGK model. The
main point of this scheme is to decompose the distribution function
on the boundary node into equilibrium and non-equilibrium parts.
Based on the mathematical analysis, we use
$-\tau_*f_i^{(1)}(t,\vec{x}_p)$ to approximate the non-equilibrium.
We tested this new boundary scheme in Couette flow and lid driven
flow. Numerical test showed that the numerical solutions of the BGK
model applied to proposed boundary scheme is very agreement
with the former paper\cite{doi:10.1146/annurev.fluid.30.1.329,1009-1963-11-4-310} and analytical
solution and also showed second order accuracy.

The proposed scheme include the different relaxation time $\tau_*$
that depend on the different viscosity between the flow and the
boundary. For this reason, we can adjust the scheme according to
viscosity between boundary and inner flow. By the equation
(\ref{boundary_solution_expand}), we can theoretically develop this
scheme to higher order accuracy. The scheme is very simply and doesn't need any more node. In one word, comparing with
the original extrapolation schemes, we have proposed a new method for the implementation of
boundary conditions for BGK model which are shown to be of second
order accuracy, and have better numerical stability.

%
%\bibliographystyle{plain}
%\bibliography{C:/Users/P/Desktop/²©Ê¿/±ÏÒµÂÛ/Reference/MyReference}

\begin{thebibliography}{10}

\bibitem{PhysRev.94.511}
P.~L. Bhatnagar, E.~P. Gross, and M.~Krook.
\newblock A {M}odel for {C}ollision {P}rocesses in {G}ases. {I}. {S}mall
  {A}mplitude {P}rocesses in {C}harged and {N}eutral {O}ne-{C}omponent
  {S}ystems.
\newblock {\em Phys. Rev.}, 94:511--525, May 1954.

\bibitem{bouzidi:3452}
M'hamed Bouzidi, Mouaouia Firdaouss, and Pierre Lallemand.
\newblock Momentum transfer of a {B}oltzmann-lattice fluid with boundaries.
\newblock {\em Physics of Fluids}, 13(11):3452--3459, 2001.

\bibitem{cercignani1988boltzmann}
C.~Cercignani.
\newblock {\em The {B}oltzmann {E}quation and {I}ts {A}pplications}.
\newblock Applied mathematical sciences. Springer-Verlag, 1988.

\bibitem{doi:10.1146/annurev.fluid.30.1.329}
Shiyi Chen and Gary~D. Doolen.
\newblock Lattice boltzmann method for fluid flows.
\newblock {\em Annual Review of Fluid Mechanics}, 30(1):329--364, 1998.

\bibitem{BOOKLi2}
D.Q.Li and T.H.Qing.
\newblock {\em Physics and Partial Differential Equations}.
\newblock Higher Education Publishing House, Beijing, 2005.

\bibitem{Filippova1998219}
Olga Filippova and Dieter H?nel.
\newblock Grid refinement for lattice-bgk models.
\newblock {\em Journal of Computational Physics}, 147(1):219 -- 228, 1998.

\bibitem{PhysRevE.67.066709}
Zhaoli Guo and T.~S. Zhao.
\newblock Explicit finite-difference lattice boltzmann method for curvilinear
  coordinates.
\newblock {\em Phys. Rev. E}, 67:066709, Jun 2003.

\bibitem{Mei1999307}
Renwei Mei, Li-Shi Luo, and Wei Shyy.
\newblock An accurate curved boundary treatment in the lattice boltzmann
  method.
\newblock {\em Journal of Computational Physics}, 155(2):307 -- 330, 1999.

\bibitem{Mei1998426}
Renwei Mei and Wei Shyy.
\newblock On the finite difference-based lattice boltzmann method in
  curvilinear coordinates.
\newblock {\em Journal of Computational Physics}, 143(2):426 -- 448, 1998.

\bibitem{Qi2001107}
Dewei Qi.
\newblock Simulations of fluidization of cylindrical multiparticles in a
  three-dimensional space.
\newblock {\em International Journal of Multiphase Flow}, 27(1):107 -- 118,
  2001.

\bibitem{0295-5075-17-6-001}
Y.~H. Qian, D.~D'Humières, and P.~Lallemand.
\newblock Lattice {BGK} {M}odels for {N}avier-{S}tokes equation.
\newblock {\em EPL (Europhysics Letters)}, 17(6):479, 1992.

\bibitem{BOOKHe1}
Y.L.He, Y.Wang, and Q.Li.
\newblock {\em Lattic {B}oltzmann method theory and applications}.
\newblock Science Publishing House, Beijing, 2008.

\bibitem{1009-1963-11-4-310}
Guo Zhao-Li, Zheng Chu-Guang, and Shi Bao-Chang.
\newblock Non-equilibrium extrapolation method for velocity and pressure
  boundary conditions in the lattice boltzmann method.
\newblock {\em Chinese Physics}, 11(4):366, 2002.

\bibitem{BOOKGuo1}
Z.L.Guo and C.G.Zheng.
\newblock {\em Theory and {A}pplications of {L}attice {B}oltzmann {M}ethod}.
\newblock Science Publishing House, Beijing, 2008.

\end{thebibliography}

\end{document}